 \documentstyle[graphicx,amscd,amssymb,verbatim,12pt,righttag,color]{amsart}
 \newtheorem{theorem}{Theorem}[section]
 \newtheorem{Def}[theorem]{Definition}
 \newtheorem{Prop}[theorem]{Proposition}
 \newtheorem{Lem}[theorem]{Lemma}
 
 \newtheorem{Rem}[theorem]{Remark}
 \newtheorem{Example}[theorem]{Example}
 \newtheorem{Coj}[theorem]{Conjecture}

 \numberwithin{equation}{section}
 \renewcommand{\rm}{\normalshape}
 
\begin{document}

\title  {Non-spectral fractal measures with Fourier frames}

\author{Chun-Kit Lai}

\address{Department of Mathematics, San Francisco State University,
1600 Holloway Avenue, San Francisco, CA 94132.}

 \email{cklai@@sfsu.edu}

\author{Yang Wang}

\address{Department of Mathematics, Hong Kong University of Science and Technology, Hong Kong}

 \email{yangwang@@ust.hk}



\begin{abstract}
We generalize the compatible tower condition  given by Strichartz to  the almost-Parseval-frame tower and show that non-trivial examples of almost-Parseval-frame tower exist. By doing so, we demonstrate  the first singular fractal measure which has only finitely many mutually orthogonal exponentials (and hence it does not admit any exponential orthonormal bases), but it still admits Fourier frames.
\end{abstract}

\maketitle

\tableofcontents
\section{Introduction}

Let $\mu$ be a compactly supported Borel probability measure on ${\Bbb R}^d$. We say that $\mu$ is a {\it frame spectral measure} if there exists a collection of exponential functions $\{e^{2\pi i \langle\lambda,x\rangle}\}_{\lambda\in\Lambda}$ such that there exists $0< A\le B<\infty$ with
$$
A\|f\|_{2}^2\leq \sum_{\lambda\in\Lambda}\left|\int f(x)e^{-2\pi i \langle\lambda,x\rangle}d\mu(x)\right|^2\leq B\|f\|_2^2, \ \forall \ f\in L^2(\mu).
$$
Whenever such $\Lambda$ exists, $\{e^{2\pi i \langle\lambda,x\rangle}\}_{\lambda\in\Lambda}$ is called a {\it Fourier frame} for $L^2(\mu)$ and $\Lambda$ is a {\it frame spectrum} for $\mu$. When $\mu$ admits an exponential orthonormal basis, we say that $\mu$ is a {\it spectral measure} and the corresponding frequency set $\Lambda$ is called a {\it spectrum} for $\mu$.

\medskip

Frames on a general Hilbert space was introduced by Duffin and Schaeffer \cite{[DS]} and it is now a fundamental building block in applied harmonic analysis. People regard frames as ``overcomplete basis" and because of its redundancy, it makes the reconstruction  more robust to errors in data and it is now widely used in signal transmission and reconstruction.  Reader may refer to \cite{[Chr]} for the background of general frame theory and \cite{[CK]} for some recent active topics.

\medskip

  One of the major hard problems in frame theory perhaps is  constructing Fourier frames or exponential orthonomal bases in different measure space $L^2(\mu)$, particularly when $\mu$ is a singular measure without any atoms or people termed it as a ``fractal measure" as the support is a fractal set.  These constructions allow Fourier analysis to work on fractal space. This problem dates back to the time of Fuglede \cite{[Fu]} who initiated the study and proposed the well-known spectral set conjecture. Although the conjecture was proved to be false by Tao \cite{[T]}, the conjecture has been extended into different facet and related questions are still being studied \cite{[W],[IKT1],[IKT2],[K],[KN]}.  Another major advance in which fractals were involved was due to Jorgensen and Pedersen \cite{[JP1]}, who discovered that the standard one-third Cantor measure is not a spectral measure, while the standard one-fourth Cantor measure is. Following the discovery, more fractal measures were found to be spectral by many others \cite{[St1],[LaW1],[DJ1]}. Many unexpected properties of the Fourier bases were discovered \cite{[St2],[DHS],[DaHL]}. While Fourier analysis appears to work perfectly on fractal spectral measures,  for the measures which are non-spectral, it is natural to ask the following question.

\medskip

{\bf (Q):} Can a non-spectral fractal measure still admit some Fourier frames?

\medskip

This question was possibly first proposed by Strichartz \cite[p.212]{[St1]}. In particular, there has been discussions  asking whether specifically the one-third Cantor measure can be frame spectral. Although we are unable to settle the case of the one-third Cantor measure, the main purpose of this paper is to answer positively {\bf (Q)} with explicit examples. (see Theorem \ref{th0.4}).

\medskip

{\bf (Q)} in its absolutely continuous counterpart is trivial since every bounded Borel set $\Omega$ with positive finite Lebesgue measure can be covered by a square. The orthonormal basis on the square naturally induces a tight frame on $\Omega$. If $\mu = g(x)dx$ is a general absolutely continuous measure, a complete characterization on the density for $\mu$ to be frame spectral was also given by the first named author \cite{[Lai]}. Such question becomes much more difficult if $\Omega$ is unbounded but still of finite measure as there cannot be any ad hoc ``square-covering" argument to construct the Fourier frames. Despite the difficulty, it was recently solved to be positive by Nitzan et al \cite{[NOU]}  who used the recent solution of the celebrated Kadison-Singer conjecture \cite{[MSS]}.

\medskip

Fractal measures are mostly supported on Lebesgue measure zero set, the situation is similar to unbounded sets of finite measures. However, it is even more complicated because if any such frame spectrum exists, there cannot be any Beurling density \cite{[DHSW],[DaHL]}. This prevents any weak convergence argument of discrete sets from happening.  Furthermore, some fractal measures are known  not to admit any Fourier frames if the measures are  non-uniform on the support \cite{[DL1]}.
Intensive researches on this question \cite{[DHSW],[DHW1],[DHW2],[DL1],[DL2],[HLL]} has been going on and one major advance was obtained recently in \cite{[DL2]}. Dutkay and Lai introduced the {\it almost-Parseval-frame condition} for the self-similar measure and proved that if such condition is satisfied, the self-similar measure admits a Fourier frame. We slightly modify the definition as below to suit the need in the paper.

\medskip

\begin{Def}
{\rm Let $\epsilon_j$ be such that $0\le\epsilon_j<1$ and $\sum_{j=1}^{\infty}\epsilon_j<\infty$. We say that $\{(N_j,B_j)\}$ is an {\it almost-Parseval-frame tower} associated to $\{\epsilon_j\}$ if}
\begin{enumerate}
\item {\rm $N_j$ are integers and $N_j\geq2$ for all $j$};
\item {\rm $B_j\subset \{0,1,...,N_j-1\}$ and $ 0 \in B_j$ for all $j$;}
\item {\rm Let $M_j := \#B_j$. There exists $L_j\subset {\Bbb Z}$ (with $0\in L_j$) such that for all $j$,}
 \begin{equation}\label{ALmost}
(1-\epsilon_j)^2\sum_{b\in B_j}|w_b|^2\leq \sum_{\lambda\in L_j}\left|\frac{1}{\sqrt{M_j}}\sum_{b\in B_j}w_be^{-2\pi i b \lambda/N_j}\right|^2\leq (1+\epsilon_j)^2\sum_{b\in B_j}|w_b|^2
\end{equation}
{\rm for all ${\bf w} = (w_b)_{b\in{B_j}}\in{\mathbb C}^{M_j}$.}
{\rm Letting the matrix ${\mathcal F}_j = \frac{1}{\sqrt{M_j}}\left[e^{2\pi i b \lambda/N_j}\right]_{\lambda\in L_j,b\in B_j}$ and $\|\cdot\|$ the standard Euclidean norm, (\ref{ALmost}) is equivalent to}
\begin{equation}\label{ALmost1}
(1-\epsilon_j)\|{\bf w}\|\leq \|{\mathcal F}_j{\bf w}\|\leq (1+\epsilon_j)\|{\bf w}\|
\end{equation}
{\rm for all  ${\bf w}\in {\Bbb C}^{M_j}$.}
 \end{enumerate}
{\rm  Whenever  $\{L_j\}_{j\in{\mathbb Z}}$ exists, we call $\{L_j\}_{j\in{\mathbb Z}}$ a {\it pre-spectrum} for the almost-Parseval-frame tower.
We define the following measures associated to an almost-Parseval-frame tower.}
$$
\nu_j = \frac{1}{M_{j}}\sum_{b\in B_{j}}\delta_{b/N_1N_2...N_j}
$$
{\rm (we denote by $\delta_a$ the Dirac measure supported on $a$) and}
\begin{equation}\label{measure}
\mu = \nu_1\ast\nu_2\ast....
\end{equation}

\end{Def}

\medskip

Roughly speaking, almost-Parseval-frame towers ensure every finite level approximated measure of the fractal is a frame spectral measure. Moreover, the frame bounds remain finite under iterations. Once all finite level has a frame with uniform frame bound, we take the weak limit under a mild condition so that fractal singular measure is also frame-spectral. 

\medskip

When all $\epsilon_j=0$, the condition is equivalent to the {\it compatible tower condition} introduced by Strichartz \cite{[St1]}. This is known to be the key condition to construct fractal spectral measures.  The measures in (\ref{measure}) are also known as {\it Moran-type measures}. These measures have been widely used in multifractal analysis \cite{[FL],[FLW]}, harmonic analysis, particularly the construction of Salem sets \cite{[LP1],[LP2]}. Some Moran-type measures were found to be spectral \cite{[AH]} and it was found later spectral Moran measures have a far reaching consequence in understanding the spectral set conjecture \cite{[GL]} and Hausdorff dimension of the support of the spectral measures \cite{[DaS]}. We note that Moran-type measures covers self-similar measures because if there exists an integer $N\ge2$ and a set $B\subset\{0,1,...,N-1\}$ such that
$$
N_j = N^{n_j}, B_j= B+NB+...+N^{n_j-1}B,
$$
 then the associated measure is the {\it self-similar measure}. In particular if $N=3$ and $B = \{0,2\}$, $\mu$ is the standard one-third Cantor measure. In such situation, the almost-Parseval-frame tower is called {\it self-similar}.

\medskip

In \cite{[DL2]}, it was proved if the almost-Parseval-frame tower is {\it self-similar}, then the self-similar measure induced will admit an Fourier frame. However, there was no example of such towers for which $\epsilon_j>0$. In this paper, we relax the self-similar restriction and  produce the first example almost-Parseval-frame tower whose $\epsilon_j>0$. We prove

\begin{theorem}\label{th0.1} Let $N_j$ and $M_j$ be positive integers satisfying
\begin{equation}\label{eq0.2}
N_j = M_jK_j+\alpha_j
\end{equation}
for some integer $K_j$ and $0\leq\alpha_j<M_j$ with
\begin{equation}\label{eq0.3}
\sum_{j=1}^{\infty}\frac{\alpha_j\sqrt{M_j}}{K_j}<\infty.
\end{equation}
Define
\begin{equation}\label{eq0.4}
B_j = \{0, K_j,...,(M_j-1)K_j\}, \ L_j = \{0,1,...,M_j-1\}.
\end{equation}
Then   $(N_{j}, B_{j})$ forms an almost-Parseval-frame tower associated with
$$
\epsilon_j=\frac{2\pi\alpha_j\sqrt{M_j}}{{K_j}}
$$ and its pre-spectrum is $\{L_j\}$. 
\end{theorem}

\medskip

We then extend the result of \cite{[DL2]} to general almost-Parseval-frame tower. For the measure $\mu$ defined in (\ref{measure}), we let
$$
\mu_{n} = \nu_1\ast...\ast\nu_{n}, \ \mu_{>n} = \nu_{n+1}\ast\nu_{n+2}\ast....
$$
so that $\mu = \mu_n\ast\mu_{>n}$. Define also the Fourier transform of a measure $\mu$ in an usual way.
$$
\widehat{\mu}(\xi) = \int e^{-2\pi i \xi x}d\mu(x).
$$

\begin{theorem}\label{th0.2}
(a) Suppose that $\{(N_j,B_j)\}$ is an almost-Parseval-frame tower associated with $\{\epsilon_j\}$ and $\{L_j\}_{j=1}^{\infty}$. Let $$
{\bf L}_n = L_1+N_1L_2+...+(N_1...N_{n-1})L_n, \ \mbox{and} \ \Lambda = \bigcup_{n=1}^{\infty}{\bf L}_n.
$$
If
$$
\delta(\Lambda) := \inf_{n}\inf_{\lambda\in{\bf L}_n}|\widehat{\mu_{>n}}(\lambda)|^2>0,
$$
 then the measure $\mu$ in (\ref{measure}) admits a Fourier frame with frame spectrum $\Lambda$.

\bigskip

(b) For the almost-Parseval-frame tower  constructed in Theorem \ref{th0.1}, the associated $\Lambda$ satisfies $\delta(\Lambda)>0 $ and hence the measure $\mu$ is a frame spectral measure.
\end{theorem}

In the end, using the theorems above,  we construct the first kind of the following examples:

\begin{theorem}\label{th0.4}
There exists non-spectral fractal measure with only finitely many orthogonal exponentials, but it still admits Fourier frames.
\end{theorem}

We organize our paper as follows: In section 2, we prove the existence of the almost-Parseval-frame tower and prove Theorem \ref{th0.1}. In Section 3, we construct the Fourier frame given the tower and prove Theorem \ref{th0.2}. In Section 4, we construct the non-spectral measures with Fourier frames. In the appendix, we study the Hausdorff dimension of the support.

\medskip

\section{Existence of Almost-Parseval-frame tower}

Let $A$ be an $n\times n$ matrix. We define the {\it operator norm} of $A$ to be
$$
\|A\| = \max_{\|{\bf x}\|=1}\|A{\bf x}\|
$$
and the {\it Frobenius norm} of $A$ to be
$$
\|A\|_F = \sqrt{\sum_{i=1}^n\sum_{j=1}^n|a_{i,j}|^2}
$$
It follows easily from Cauchy-Schwarz inequality that $\|A\|_2\leq \|A\|_{F}$.

\medskip

For $N_j$ and $M_j$ satisfying (\ref{eq0.2}) and (\ref{eq0.3}) and for $B_j$ and $L_j$ defined in (\ref{eq0.4}),
We let
$$
{\mathcal F}_j= \frac{1}{\sqrt{M_j}}\left[e^{2\pi i b \lambda/N_j}\right]_{\lambda\in L_j,b\in B_j}, \ {\mathcal H}_j= \frac{1}{\sqrt{M_j}}\left[e^{2\pi i b \lambda/M_jK_j}\right]_{\lambda\in L_j,b\in B_j}.
$$

\medskip

\begin{Lem}\label{lem1.1}
${\mathcal H}_n$ is a unitary matrix. i.e. $\|{\mathcal H}_n{\bf x}\| = \|{\bf x}\|$.
\end{Lem}

\begin{pf}
Let $b = m\in L_j$ and $\lambda = nK_j\in B_j$, for $m,n = 0,1,...,M_j-1$. It follows directly that  $e^{2\pi i b \lambda/MK_n} = e^{2\pi i mn/M_j}$. Hence,
$$
{\mathcal H}_j = \frac{1}{\sqrt{M_j}}\left[e^{2\pi i mn/M_j}\right]_{m,n = 0,...,M_j-1},
$$
which is the standard Fourier matrix of order $M_j$. Thus, ${\mathcal H}_j$ is unitary.
\end{pf}

\medskip

\noindent{\it Proof of Theorem \ref{th0.1}.} We first show that for any $j>0$,
\begin{equation}\label{eq2.0}
\|{\mathcal F}_j-{\mathcal H}_j\|\leq\frac{2\pi \alpha_j\sqrt{M_j}}{K_j}.
\end{equation}
To see this, We note that
\begin{equation}\label{eq2.1}
\|{\mathcal F}_j-{\mathcal H}_j\|^2\leq \|{\mathcal F}_j-{\mathcal H}_j\|_{F}^2 = \frac{1}{M_j}\sum_{b\in B_j}\sum_{\lambda\in L_j}\left|e^{2\pi i b\lambda/N_j}-e^{2\pi i b\lambda/M_jK_j}\right|^2.
\end{equation}
We now estimate the difference of the exponentials inside the summation. Recall that for any $\theta_1,\theta_2$,
$$
|e^{i\theta_1}-e^{i\theta_2}| = |e^{i(\theta_1-\theta_2)}-1|\leq |\theta_1-\theta_2|.
$$
This implies that
$$
\begin{aligned}
\left|e^{2\pi i b\lambda/N_j}-e^{2\pi i b\lambda/M_jK_j}\right|^2\leq& \left|\frac{2\pi b\lambda}{N_j}-\frac{2\pi b\lambda}{M_jK_j}\right|^2\nonumber\\
=& 4\pi^2\frac{b^2\lambda^2 \alpha_j^2}{M_j^2K_j^2N_j^2}  \ \ \  \ \mbox{(by $N_j = M_jK_j+\alpha_j$)}\nonumber\\
\leq &4\pi^2\frac{M_j^2\alpha_j^2}{N_j^2}  \ \ \ \ \mbox{(by $b\leq M_jK_j$ and $\lambda\leq M_j$)}\\
\end{aligned}
$$
Hence, from (\ref{eq2.1}),
\begin{eqnarray}\label{eq2.2}
\|{\mathcal F}_j-{\mathcal H}_j\|^2&\leq& \frac{1}{M_j}\sum_{b\in B_j}\sum_{\lambda\in L_j}4\pi^2\frac{M_j^2\alpha_j^2}{N_j^2} \nonumber\\
 &=& 4\pi^2\frac{M_j^3\alpha_j^2}{N_j^2} \nonumber \\
 &= & 4\pi^2\frac{M_j\alpha_j^2}{\left( K_j+\alpha_j/M_j\right)^2}
\end{eqnarray}
As $\alpha_j\ge 0$, $\|{\mathcal F}_j-{\mathcal H}_j\|^2\leq  4\pi^2{\alpha_j^2 M_j}/{K_j^2}$ and thus  (\ref{eq2.0}) follows by taking square root.

\medskip

 We now show that $\{(N_{j}, B_{j})\}$ forms an almost-Parseval-frame tower with pre-spectrum $L_j$. The first two conditions for the almost-Parseval-frame tower are clearly satisfied. To see the last condition, we recall that $\epsilon_j = 2\pi \sqrt{M_j}\alpha_j/{K_j} $.  From the triangle inequality and (\ref{eq2.0}), we have
$$
\begin{aligned}
\|{\mathcal F}_j{\bf w}\|\leq&\|{\mathcal H}_j{\bf w}\| + \|{\mathcal F}_j-{\mathcal H}_j\|\|{\bf w}\|\\
\leq&  \left(1+\frac{2\pi \alpha_j\sqrt{M_j}}{K_j}\right)\|{\bf w}\| = (1+\epsilon_j)\|{\bf w}\|.\\
\end{aligned}
$$
Similarly, for the lower bound,
$$
\begin{aligned}
\|{\mathcal F}_{j}{\bf w}\|\geq &\|{\mathcal H}_{j}{\bf w}\|-\|{\mathcal F}_{j}-{\mathcal H}_{j}\|\|{\bf w}\|\\
\geq& \left(1-\frac{2\pi \alpha_j\sqrt{M_j}}{K_j}\right)\|{\bf w}\| = (1-\epsilon_j)\|{\bf w}\|.
\end{aligned}
$$
Thus, from (\ref{ALmost1}), the last condition follows and $(N_j,B_j)$ satisfies the almost-Parseval-frame condition associated with $\{\epsilon_j\}$ and $\sum_{j=1}^{\infty}\epsilon_j<\infty$ is guaranteed by (\ref{eq0.3}) in the assumption.
\qquad$\Box$

\begin{Rem}
In view of (\ref{eq2.2}), condition (\ref{eq0.3}) can be replaced by a weaker condition
$$
\sum_{j=1}^{\infty} \frac{\alpha_j\sqrt{M_j}}{ K_j+\alpha_j/M_j}<\infty.
$$
(\ref{eq0.3}) would be enough  for the convenience of our discussion. It is also worth to note that if all $\alpha_j=0$, then the the matrices ${\mathcal F}_j = {\mathcal H}_j$ are reduced to the Hadamard matrices. The associated measures are all spectral measures. see e.g. \cite{[AH]}.
\end{Rem}

\medskip

We end this section by illustrating some explicit examples of Theorem \ref{th0.1}.

\begin{Example}\label{Example2.3}
{\rm Let $p$ be an odd prime and  suppose that $N_j = p^j$. Let $M_j = 2$ for all $j$. Then it is clear that $N_j = 2K_j+1$ for some $K_j$. In this case,}
$$
\sum_{j=1}^{\infty} \frac{\alpha_j\sqrt{M_j}}{K_j} = \sum_{j=1}^{\infty} \frac{\sqrt{2}}{K_j} = \sum_{j=1}^{\infty} \frac{2\sqrt{2}}{p^j-1}<\infty.
$$
{\rm Thus $N_j = p^j$ and $B_j = \{0,K_j\}$ forms an almost-Parseval-frame tower with pre-spectrum $L_j = \{0,1\}$ for all $j$.}
\end{Example}

\medskip

\begin{Example}\label{Example2.4}
 For $0\le\beta<2$, $\gamma\ge 0$ and $N\ge 1$  such that $$
N(1-\beta/2-\gamma)>1,$$
 let $K_j,M_j,\alpha_j$ be integers $K_j\ge j^N$, $M_j \leq K_j^{\beta}$ and $\alpha_j \le  K_j^{\gamma}$. Then
$$
\sum_{j=1}^{\infty} \frac{\alpha_j\sqrt{M_j}}{K_j}\leq \sum_{j=1}^{\infty} \frac{K_j^{\gamma}K_j^{\beta/2}}{K_j} = \sum_{j=1}^{\infty} \frac{1}{K_j^{1-\gamma-\beta/2}}\le\sum_{j=1}^{\infty} \frac{1}{j^{N(1-\gamma-\beta/2)}}<\infty.
$$
Hence, $N_j = K_jM_j+\alpha_j$ and $B_j = \{0,K_j,...,(M_j-1)K_j\}$ satisfies the almost-Parseval-frame condition.
\end{Example}

\medskip

\section{Construction of Fourier frames}

In this section, we consider the almost-Parseval-frame tower defined in Section 1 and show that the measure $\mu$ defined in (\ref{measure}) is a frame spectral measure. We first recall some notations.
$$
\nu_j = \frac{1}{M_j}\sum_{b\in B_j}\delta_{b/N_1...N_j}, \ \mbox{and} \ \mu = \nu_1\ast\nu_2....
$$
We define
$$
\mu_{ n} = \nu_1\ast...\ast\nu_n, \ \mu_{>n}= \nu_{n+1}\ast\nu_{n+2}\ast...
$$
so that $\mu  =\mu_{ n}\ast\mu_{>n}$. It is also direct to see that the support of $\mu$ is the compact set
$$
K_{\mu} = \left\{\sum_{j=1}^{\infty}\frac{b_j}{N_1...N_j}: b_j\in B_j \ \mbox{for all} \ j\right\}.
$$
We also consider the first $n^{\rm th}$-partial sum in $K_{\mu}$ and denote it by
$$
{\bf B}_n = \frac{1}{N_1}B_1+\frac{1}{N_1N_2}B_2+...+\frac{1}{N_1N_2...N_n}B_n
$$
which is the support of $\mu_n$. For the $\{L_j\}_{j\in{\mathbb Z}}$ in the tower, we consider
$$
{\bf L}_n = L_1+N_1L_2+...+(N_1...N_{n-1})L_n.
$$

\medskip

\begin{Prop}\label{prop3.1}
For any $n\ge 1$, let ${\bf M}_n= \prod_{j=1}^{n}M_j$ we have
$$
\left(\prod_{j=1}^{n}(1-\epsilon_j)\right)^2\|{\bf w}\|^2\leq\sum_{\lambda\in{\bf L}_n}\left|\frac{1}{\sqrt{{\bf M}_{n}}}\sum_{{\bf b}\in{\bf B}_n}w_{\bf b}e^{-2\pi i {\bf b}\lambda}\right|^2\leq\left(\prod_{j=1}^{n}(1+\epsilon_j)\right)^2\|{\bf w}\|^2
$$
for any ${\bf w} = (w_{{\bf b}})_{{\bf b}\in{\bf B}_n}\in{\mathbb C}^{M_1...M_n}$
\end{Prop}

\begin{pf}
We prove it by mathematical induction. When $n=1$, it is the almost-Parseval condition for $(N_1,B_1)$ so the statement is true trivially. Assume now the inequality is true for $n-1$. Then we decompose ${\bf b}\in {\bf B}_n$ and $\lambda\in {\bf L}_n$ by
$$
{\bf b} = \frac{1}{N_1...N_n}b_n+{\bf b}_{n-1}, \ \lambda = \lambda_{n-1}+ N_1...N_{n-1}l_n,
$$
where $b_n\in B_{n}$, ${\bf b}_{n-1}\in{\bf B}_{n-1}$, $\lambda_{n-1}\in {\bf L}_{n-1}$ and $l_n\in L_n$. Now, we have
$$
\begin{aligned}
&\sum_{\lambda\in{\bf L}_n}\left|\sum_{{\bf b}\in{\bf B}_{n}}w_{\bf b}e^{-2\pi i {\bf b}\lambda}\right|^2\\
=&\sum_{\lambda_{n-1}\in{\bf L}_{n-1}}\sum_{l_n\in L_n}\left|\sum_{{\bf b}_{n-1}\in{\bf B}_{n-1}}\sum_{b_n\in B_n}\frac{1}{\sqrt{{\bf M}_{n}}}w_{{\bf b}b_n}e^{-2\pi i \left( \frac{1}{N_1...N_n}b_n+{\bf b}_{n-1}\right)\cdot\left(\lambda_{n-1}+ N_1...N_{n-1}l_n\right)}\right|^2.\\
\end{aligned}
$$
Note that ${\bf b}_{n-1}\cdot(N_1...N_{n-1})l_n $ is always an integer, the right hand side above can be written as
$$
\sum_{\lambda_{n-1}\in{\bf L}_{n-1}}\sum_{l_n\in L_n}\left|\sum_{b_n\in B_n}\frac{1}{\sqrt{{ M}_{n}}}\left(\sum_{{\bf b}_{n-1}\in{\bf B}_{n-1}}\frac{1}{\sqrt{{\bf M}_{n-1}}}w_{{\bf b}_{n-1}b_n}e^{-2\pi i \left( \frac{1}{N_1...N_n}b_n+{\bf b}_{n-1}\right)\cdot\lambda_{n-1}}\right)e^{-2\pi i  b_nl_n/N_n}\right|^2
$$
Using the almost-Parseval-frame condition for $(N_n,B_n)$ and also the induction hypothesis, this term
$$
\begin{aligned}
\leq&(1+\epsilon_n)^2\sum_{\lambda_{n-1}\in{\bf L}_{n-1}}\sum_{b_n\in B_n}\left|\sum_{{\bf b}_{n-1}\in{\bf B}_{n-1}}\frac{1}{\sqrt{{\bf M}_{n-1}}}w_{{\bf b}_{n-1}b_n}e^{-2\pi i \left( \frac{1}{N_1...N_n}b_n+{\bf b}_{n-1}\right)\cdot\lambda_{n-1}}\right|^2\\
=&(1+\epsilon_n)^2\sum_{b_n\in B_{n}}\sum_{\lambda_{n-1}\in{\bf L}_{n-1}}\left|\sum_{{\bf b}_{n-1}\in{\bf B}_{n-1}}\frac{1}{\sqrt{{\bf M}_{n-1}}}w_{{\bf b}_{n-1}b_n}e^{-2\pi i {\bf b}_{n-1}\cdot\lambda_{n-1}}\right|^2\\
\leq &\left(\prod_{j=1}^{n}(1+\epsilon_j)\right)^2\sum_{b_n\in B_{n}}\sum_{{\bf b}_{n-1}\in{\bf B}_{n-1}}|w_{{\bf b}_{n-1}b_n}|^2\\
=&\left(\prod_{j=1}^{n}(1+\epsilon_j)\right)^2\|{\bf w}\|^2.\\
\end{aligned}
$$
This completes the proof of the upper bound and the proof of the lower bound is analogous.
\end{pf}

\medskip

We now decompose  $K_{\mu}$ as
\begin{equation}\label{K_mu}
K_{\mu} = \bigcup_{{\bf b}\in {\bf B}_n} ({\bf b}+ K_{\mu,n}).
\end{equation}
where
$$
K_{\mu,n} = \left\{\sum_{j={n+1}}^{\infty}\frac{b_j}{N_1...N_j}: b_j\in B_j \ \mbox{for all} \ j\right\}.
$$
 Denote by $K_{\bf b} = {\bf b}+ K_{\mu,n}$ and ${\bf 1}_{K_{\bf b}}$ the characteristic function of $K_{\bf b}$. Let
$$
{\mathcal S}_n = \left\{\sum_{{\bf b}\in {\bf B}_n}w_{\bf b}{\bf 1}_{K_{\bf b}}: w_{\bf b}\in{\mathbb C}\right\}.
$$
${\mathcal S}_n$ denotes the collection of all $n^{th}$ level step functions on $K_{\mu}$. As
$$
K_{\mu,n} = \bigcup_{b\in B_{n+1}}\left(\frac{b}{N_1...N_{n+1}}+K_{\mu,n+1}\right),
$$
we have ${\mathcal S}_1\subset {\mathcal S}_2\subset....$. Let also
$$
{\mathcal S} = \bigcup_{n=1}^{\infty}{\mathcal S}_n.
$$
It is clear that ${\mathcal S}$ forms a dense set of functions in $L^2(\mu)$.
\medskip

\begin{Lem}\label{lem3.1}
Let $f = \sum_{{\bf b}\in {\bf B}_n}w_{\bf b}{\bf 1}_{K_{\bf b}}\in {\mathcal S}_n$. Then
\begin{equation}\label{eq3.1}
\int|f|^2d\mu = \frac{1}{{\bf M}_n}\sum_{{\bf b}\in {\bf B}_n}|w_{\bf b}|^2.
\end{equation}
\begin{equation}\label{eq3.2}
\int f(x)e^{-2\pi i \lambda x}d\mu(x) = \frac{1}{{\bf M}_n}\widehat{\mu_{>n}}(\lambda)\sum_{{\bf b}\in {\bf B}_n}w_{\bf b} e^{-2\pi i  {\bf b} \lambda}.
\end{equation}
Here ${\bf M}_n = M_1...M_n$.
\end{Lem}

\begin{pf}
As $K_{\bf b}$ and $K_{\bf b}'$ has either empty intersection or intersects at most one point, taking $\mu$-measure on (\ref{K_mu}), we obtain  $\mu(K_{\bf b}) = 1/{\bf M}_n$.  (\ref{eq3.1}) follows from a direct computation. For (\ref{eq3.2}), we use $\mu = \mu_n\ast\mu_{>n}$ and we have
$$
\begin{aligned}
\int f(x)e^{-2\pi i \lambda x}d\mu(x) = &\sum_{{\bf b}\in{\bf B}_n}w_{\bf b} \int {\bf 1}_{K_{\bf b}}(x)e^{-2\pi i \lambda x}d(\mu_n\ast\mu_{>n}(x)) \\
=& \sum_{{\bf b}\in{\bf B}_n}w_{\bf b} \int {\bf 1}_{{\bf b}+K_{\mu,n}}(x+y)e^{-2\pi i \lambda (x+y)}d\mu_n(x)d\mu_{>n}(y).\\
\end{aligned}
$$
Note that $\mu_{>n}$ is supported on $K_{\mu,n}$ and $K_{\bf b}$ and $K_{\bf b}'$ has either empty intersection or intersects at most one point. The above is equal to
$$
\begin{aligned}
=&\sum_{{\bf b}\in{\bf B}_n}w_{\bf b} \frac{1}{{\bf M}_n}\int {\bf 1}_{{\bf b}+K_{\mu,n}}({\bf b}+y)e^{-2\pi i \lambda ({\bf b}+y)}d\mu_{>n}(y)\\
=&\sum_{{\bf b}\in{\bf B}_n}w_{\bf b} \frac{1}{{\bf M}_n}e^{-2\pi i \lambda{\bf b}}\int e^{-2\pi i \lambda y}d\mu_{>n}(y)\\
=& \frac{1}{{\bf M}_n}\widehat{\mu_{>n}}(\lambda)\sum_{{\bf b}\in {\bf B}_n}w_{\bf b} e^{-2\pi i  {\bf b} \lambda}.\\
\end{aligned}
$$
The lemma follows.
\end{pf}

\medskip

Let \begin{equation}\label{Lambda}
\Lambda = \bigcup_{n=1}^{\infty} {\bf L}_n.
\end{equation}
As $0\in L$, the sets in the union is an increasing union. We now define the following quantity
$$
\delta_{n_0}(\Lambda) = \inf_{n\ge n_0}\inf_{\lambda\in {\bf L}_n}|\widehat{\mu_{>n}}(\lambda)|^2,
$$
for some $n_0\ge 1$. The following theorem gives a sufficient condition for $\Lambda$ to be a Fourier frame for $\mu$.

\begin{theorem}\label{th3.3}
Suppose that $(N_j,B_j)$ is an almost-Parseval-frame tower and $\mu$ be the associated measure. Let $L$ be the associated spectrum for the tower and $\Lambda$ defined (\ref{Lambda}) satisfies $\delta_{n_0}(\Lambda)>0$. Then $\mu$ admits a Fourier frame $E(\Lambda)$ with  lower and upper frame bounds respectively equal
$$
\delta(\Lambda)\left(\prod_{j=1}^{\infty}(1-\epsilon_j)\right)^2, \left(\prod_{j=1}^{\infty}(1+\epsilon_j)\right)^2.
$$
\end{theorem}

\begin{pf}
To check the Fourier frame inequality holds, it suffices to show that it is true for a dense set of functions in $L^2(\mu)$ \cite[Lemma 5.1.7]{[Chr]}, in which we will check it for step functions in ${\mathcal S}$. Moreover, since ${\mathcal S}_n$ is an increasing union of sets, we consider $f = \sum_{{\bf b}\in {\bf B}_n}w_{\bf b}{\bf 1}_{K_{\bf b}}\in {\mathcal S}_n$ with $n\ge n_0$. By Lemma \ref{lem3.1}, we have
$$
\begin{aligned}
\sum_{\lambda\in {\bf L}_{n}}\left|\int f(x)e^{-2\pi i \lambda x}d\mu(x)\right|^2 =& \sum_{\lambda\in {\bf L}_{n}}\left|\frac{1}{{\bf M}_n}\widehat{\mu_{>n}}(\lambda)\sum_{{\bf b}\in {\bf B}_n}w_{\bf b} e^{-2\pi i  {\bf b} \lambda}\right|^2\\
=& \frac{1}{{\bf M}_n}\sum_{\lambda\in {\bf L}_{n}}\left|\widehat{\mu_{>n}}(\lambda)\right|^2\left|\frac{1}{\sqrt{{\bf M}_n}}\sum_{{\bf b}\in {\bf B}_n}w_{\bf b} e^{-2\pi i  {\bf b} \lambda}\right|^2\\
\end{aligned}
$$
Note that $\delta(\Lambda)\leq \left|\widehat{\mu_{>n}}(\lambda)\right|^2\leq 1$. By Proposition \ref{prop3.1}, we have this implies that
$$
\frac{1}{{\bf M}_n}\delta(\Lambda)\left(\prod_{j=1}^{n}(1-\epsilon_j)\right)^2\|{\bf w}\|^2
\leq\sum_{\lambda\in {\bf L}_{n}}\left|\int f(x)e^{-2\pi i \lambda x}d\mu(x)\right|^2\leq  \frac{1}{{\bf M}_n} \left(\prod_{j=1}^{n}(1+\epsilon_j)\right)^2\|{\bf w}\|^2.
$$
Using Lemma \ref{lem3.1} again, we have
$$
\delta(\Lambda)\left(\prod_{j=1}^{n}(1-\epsilon_j)\right)^2\int|f|^2d\mu
\leq\sum_{\lambda\in {\bf L}_{n}}\left|\int f(x)e^{-2\pi i \lambda x}d\mu(x)\right|^2\leq  \left(\prod_{j=1}^{n}(1+\epsilon_j)\right)^2\int|f|^2d\mu.
$$
To complete the proof, we note that for all $m>n$,  $f\in {\mathcal S}_n\subset{\mathcal S}_m$, the inequality can also be written as
$$
\delta(\Lambda)\left(\prod_{j=1}^{m}(1-\epsilon_j)\right)^2\int|f|^2d\mu
\leq\sum_{\lambda\in {\bf L}_{m}}\left|\int f(x)e^{-2\pi i \lambda x}d\mu(x)\right|^2\leq  \left(\prod_{j=1}^{m}(1+\epsilon_j)\right)^2\int|f|^2d\mu,
$$
for all $f\in {\mathcal S}_n$. Taking $m$ to infinity, we have
$$
\delta(\Lambda)\left(\prod_{j=1}^{\infty}(1-\epsilon_j)\right)^2\int|f|^2d\mu
\leq\sum_{\lambda\in \Lambda}\left|\int f(x)e^{-2\pi i \lambda x}d\mu(x)\right|^2\leq  \left(\prod_{j=1}^{\infty}(1+\epsilon_j)\right)^2\int|f|^2d\mu.
$$
Note that the frame bounds are finite since $\sum_{j=1}^{\infty}\epsilon_j<\infty$. This shows the frame inequality for any $f\in{\mathcal S}$. Hence, $E(\Lambda)$ is a Fourier frame for $L^2(\mu)$.

\end{pf}

\medskip

We now show that the almost-Parseval-frame tower constructed in Theorem \ref{th0.1} satisfies $\delta(\Lambda)>0$. Recall that
$$
 N_j= M_jK_j+\alpha_j,
$$
with $B_{j} = \{0,K_j,...,(M_j-1)K_j\}$ and $L_j= \{0,1,...,M_j-1\}$. The associated measure is given by
$$
\mu = \nu_1\ast\nu_2\ast..., \ \mbox{and} \ \nu_j = \frac{1}{M_j}\sum_{b\in B_{j}}\delta_{j/N_1...N_j}.
$$
The Fourier transform is given by
$$
\widehat{\mu}(\xi) = \prod_{j=1}^{\infty}\widehat{\nu_j}(\xi) =\prod_{j=1}^{\infty}\left[\frac{1}{M_j}\sum_{k=0}^{M_j-1}e^{-2\pi i kK_{j}\xi/N_1...N_j} \right]
$$
It follows directly from summation of geometric series that
$$
\begin{aligned}
\widehat{\nu_j}(\xi)=& \left\{
     \begin{array}{ll}
       \frac{1}{M_j}e^{\pi i c_j(M_j-1)\xi}\frac{\sin \pi c_jM_j \xi}{\sin\pi c_j\xi}
, & \hbox{if $\xi\not\in \frac{1}{c_j}{\mathbb Z}$;} \\
       1, & \hbox{if $\xi\in \frac{1}{c_j}{\mathbb Z}$.}
     \end{array}
   \right.
\end{aligned}
$$
where $c_j = K_{j}/N_1...N_j$. We prove that

\medskip

\begin{Prop}\label{prop3.2}
With all the notation above, there exists $k_0$ such that $\Lambda = \bigcup_{k=1}^{\infty} {\bf L}_k$ satisfies
$$
\delta_{k_0}(\Lambda) = \inf_{k\ge k_0}\inf_{\lambda\in{\bf L}_k} |\widehat{\mu_{>k}}(\lambda)|^2>0
$$
where ${\bf L}_k = L_1+N_1L_2+...+N_1...N_{k-1}L_k$.
\end{Prop}

\begin{pf}
We note that
\begin{equation}\label{eq3.2+}
 |\widehat{\mu_{>k}}(\lambda)|^2 = \prod_{j=1}^{\infty}|\widehat{\nu_{k+j}}(\lambda)|^2 =  \prod_{j=1}^{\infty}\left|\frac{1}{M_{k+j}}\sum_{\ell=0}^{M_{k+j}-1}e^{-2\pi i \ell K_{k+j}\lambda/(N_1...N_kN_{k+1}...N_{k+j})}\right|^2.
\end{equation}
For any $\lambda\in{\bf L}_k$ for which the terms $|\widehat{\nu_{k+j}}(\lambda)|^2<1$, we have
$$
\left|\frac{1}{M_{k+j}}\sum_{k=0}^{M_j-1}e^{-2\pi i kK_{j}\lambda/(N_1...N_kN_{k+1}...N_{k+j})}\right|^2= \left|\frac{1}{M_{k+j}}\frac{\sin \pi c_{k+j}M_{k+j} \lambda}{\sin\pi c_{k+j}\lambda}\right|^2.
$$
Using the elementary estimate $\sin x\leq x$ and $\sin x\geq x-x^{3}/3!$, we have
$$
\begin{aligned}
\left|\frac{1}{M_{k+j}}\frac{\sin \pi c_{k+j}M_{k+j} \lambda}{\sin\pi c_{k+j}\lambda}\right|^2\geq &\left|\frac{\sin (\pi c_{k+j} M_{k+j}\lambda)}{\pi c_{k+j} M_{k+j}\lambda}\right|^2= \left(1-\frac{(\pi c_{k+j} M_{k+j}\lambda)^2}{3!}\right)^2\\
\end{aligned}
$$
Recall that $c_{k+j} = K_{{k+j}}/N_1...N_{k+j}$, we have
$$
\left(1-\frac{(\pi c_{k+j} M_{k+j}\lambda)^2}{3!}\right)^2 = \left(1-\frac{\pi^2}{6N_{k+1}^2...N_{k+j-1}^2}\cdot\left(\frac{K_{{k+j}}M_{k+j}}{N_{k+j}}\right)^2 \cdot \left(\frac{\lambda}{N_1...N_{k}}\right)^2\right)^2
$$
\begin{equation}\label{eq3.3}
\end{equation}
We need to ensure all the terms inside the outermost square are positive and their product is strictly positive. For $\lambda\in {\bf L}_k$, we write
$$
\lambda = \ell_1+N_1\ell_2+....+(N_1...N_{k-1})\ell_k, \ \mbox{for some} \  \ell_i\in L_i.
$$
From $N_i = M_iK_i+\alpha_i$, we have $\ell_i\leq M_i-1<N_i$  and thus
$$
\begin{aligned}
\frac{\lambda}{N_1...N_k} =&  \frac{\ell_k}{N_k}+\frac{\ell_{k-1}}{N_kN_{k-1}}+...+\frac{\ell_1}{N_k...N_1}\\
 <& \frac{M_k}{N_k}+\frac{M_{k-1}}{N_kN_{k-1}}+...+\frac{M_1}{N_k...N_1}\\
\le& \frac{1}{K_k}+\frac{1}{N_kK_{k-1}}+...+\frac{1}{N_k...N_2K_1}\\
\le& \frac{1}{K_k}\left(1+\frac{1}{M_k K_{k-1}}+\frac{1}{M_kN_{k-1}K_{k-2}}+...+\frac{1}{M_kN_{k-1}N_{k-2}...N_{2}K_1}\right)\\
\le& \frac{2}{K_{k}}. (\mbox{since all} \ M_j,N_j\ge 2)
\end{aligned}
$$
For the term $j>1$ in (\ref{eq3.3}),
we use
$$
\frac{\lambda}{N_1...N_{k}}<1, \mbox{and} \ \frac{K_{{k+j}}M_{k+j}}{N_{k+j}}\leq 1.
$$
With $N_j\ge 2$ for all $j$, we have
$$
(\ref{eq3.3}) \geq \left(1-\frac{\pi^2}{6\cdot2^{2(j-1)}}\right)^2, \ \mbox{for} \ j>1.
$$
If $j=1$,
$$
\left(1-\frac{\pi^2}{6}\cdot\left(\frac{K_{{k+j}}M_{k+j}}{N_{k+j}}\right)^2 \cdot \left(\frac{\lambda}{N_1...N_{k}}\right)^2\right)^2\ge \left(1-\frac{\pi^2}{6} \cdot \left(\frac{2}{K_k}\right)^2\right)^2\ge \left(1-\frac{2\pi^2}{3K_k^2} \right)^2
$$
Note that our assumption that $\sum_{k=1}^{\infty}\frac{\alpha_k\sqrt{M_k}}{K_k}<\infty$ implies that $K_k$ tends to infinity. Hence, there exists $k_0$ such that for all $k\ge k_0$, $K_k\ge 3$. This ensure the term insider the square is greater than or equal to $\delta: = 1-2\pi^2/27>0$. Putting all the inequality back to (\ref{eq3.2+}), we obtain
$$
 |\widehat{\mu_{>k}}(\lambda)|^2\geq \delta^2\cdot\prod_{j=2}^{\infty}\left(1-\frac{\pi^2}{6\cdot2^{2(j-1)}}\right)^2 : = c_0.
$$
Hence, $\delta(\Lambda)\geq c_0$. As $\sum_{j=2}^{\infty}{\pi^2}/{(6\cdot2^{2(j-1)})}<\infty$ and  ${\pi^2}/{(6\cdot2^{2(j-1)})}<1$ for all $j>1$, $c_0>0$ and this completes the proof.

\end{pf}

\medskip

\noindent{\it Proof of Theorem \ref{th0.2}.} (a) follows directly from Theorem \ref{th3.3}. For (b), Proposition \ref{prop3.2} implies that $\delta(\Lambda)>0$ and hence the measure $\mu$ is frame spectral by Theorem \ref{th3.3}.\qquad$\Box$

\medskip

\section{Non-spectral measures}
In this section, we will see the measures defined by the almost-Parseval-frame tower in Theorem \ref{th0.1} is in general not spectral.
For a given probability measure $\mu$, we let
$$
Z(\widehat{\mu}) = \{\xi\in{\Bbb R}: \widehat{\mu}(\xi)=0\}
$$
be its zero set of $\widehat{\mu}$. We recall that the collection of the exponentials $\{e^{2\pi i \lambda x}: \lambda\in\Lambda\}$ is a {\it mutually orthogonal set} if the exponential functions are mutually orthogonal in $L^2(\mu)$.
 In order to show $\mu$ cannot be a spectral measure, we need the following simple observation.

\begin{Lem}\label{lem1.3}
If $\mu$ is a spectral measure whose support is an infinite set, then any mutually orthogonal set $\Lambda$ must be of infinite cardinality and satisfies $\Lambda-\Lambda\subset Z(\widehat{\mu})\cup\{0\}$.
\end{Lem}

\medskip

\begin{pf}
If $\mu$ is a spectral measure whose support is an infinite set, then $L^2(\mu)$ is of infinite dimension as a vector space, so any mutually orthogonal sets must be infinite in cardinality. For mutually orthogonality to hold, we need  for any $\lambda\ne\lambda'\in\Lambda$,
$$
0=\int e^{2\pi i (\lambda-\lambda')}xd\mu(x) =\widehat{\mu}(\lambda-\lambda')
$$
Hence, $\Lambda-\Lambda\subset Z(\widehat{\mu})\cup\{0\}$ follows.
\end{pf}

\medskip

Focusing on the tower we constructed in Theorem \ref{th0.1},
$$
N_j = K_{j}M_j+\alpha_j.
$$
and $ B_{j} = \{0, K_{j},...,(M_j-1)K_{j}\}, \ L = \{0,1,...,M_j-1\},$  we have

\begin{Lem}\label{lem1.4}
$$Z(\widehat{\mu}) = \bigcup_{j=1}^{\infty}Z(\widehat{\nu}_j) = \bigcup_{j=1}^{\infty}\left[\frac{N_1...N_j}{K_{j}M_j}({\mathbb Z}\setminus M_j{\Bbb Z})\right].
$$
\end{Lem}

\begin{pf}
 we can compute directly the zero set of the Fourier transform of $\nu_j$ as
$$
\begin{aligned}
\widehat{\nu_j}(\xi) =&\frac{1}{M_j}\sum_{k=0}^{M_j-1}e^{-2\pi i kK_{j}\xi/N_1...N_j} \\
=& \left\{
     \begin{array}{ll}
       \frac{1}{M_j}e^{\pi i c_j(M_j-1)\xi}\frac{\sin \pi c_jM_j \xi}{\sin\pi c_j\xi}
, & \hbox{if $\xi\not\in \frac{1}{c_j}{\mathbb Z}$;} \\
       1, & \hbox{if $\xi\in \frac{1}{c_j}{\mathbb Z}$.}
     \end{array}
   \right.
\end{aligned}
$$
where $c_j = K_{j}/N_1....N_{j}$. It follows directly that
$$
Z(\widehat{\nu}_j) = \frac{1}{c_jM_j}({\mathbb Z}\setminus M_j{\Bbb Z})= \frac{N_1...N_j}{K_{j}M_j}({\mathbb Z}\setminus M_j{\Bbb Z})
$$
so that
$$
Z(\widehat{\mu}) = \bigcup_{j=1}^{\infty}Z(\widehat{\nu}_j) = \bigcup_{j=1}^{\infty}\left[\frac{N_1...N_j}{K_{j}M_j}({\mathbb Z}\setminus M_j{\Bbb Z})\right].
$$
\end{pf}

\medskip
It is natural to conjecture that

\begin{Coj}
Suppose that $(N_j,B_j)$ are the almost-Parseval-frame tower defined in Theorem \ref{th0.1} and the associated measure $\mu$ is spectral. Then all $\alpha_j =0$.
\end{Coj}

However, this will let us into rather involved number theoretic and combinatoric questions. To serve the purpose of this paper,  the following proposition shows that under simple conditions on $M_j$, $K_{j}$ and $\alpha_j$, the measure $\mu$ cannot be spectral.

\begin{Prop}\label{prop1.2}
Suppose that $\Lambda$ is a mutually orthogonal set for $\mu$ defined in (\ref{measure}) with $N_j = K_jM_j+1$ $(\alpha_j=1)$ and $B_{j} = \{0,K_j,...,(M_j-1)K_j\}, L_{j} = \{0,1,...,M_j-1\}$ satisfying

\medskip

(1) all $M =2$ and

\medskip

(2) all $K_j$ are odd; 

\medskip
 Then
$$
\#\Lambda\leq 2.
$$
\end{Prop}

\begin{pf}
Suppose that there exists mutually orthogonal set $\Lambda$ with cardinality  greater than $2$. We can find distinct $\lambda_1,\lambda_2,\lambda_3$ such that $\lambda_1-\lambda_2,\lambda_3-\lambda_2,\lambda_1-\lambda_3\in Z(\widehat{\mu})$. Hence, we can write
$$
\lambda_1-\lambda_2 = \frac{(N_1...N_j)}{K_{j}M_j}(r_j+ M_jq_j), \ \lambda_3-\lambda_2 = \frac{(N_1...N_k)}{K_{k}M_k}(r_k+ M_kq_k),
$$
$$ \lambda_1-\lambda_3 = \frac{(N_1...N_{\ell})}{K_{{\ell}}M_{\ell}}(r_\ell+ M_{\ell}q_\ell)
$$
 where $0<r_n<M_n$ for $n=j,k,\ell$ and $q_j,q_k,q_{\ell}$ are integers. Denote by $$
{\bf N}_n = N_1...N_n.$$
 As $(\lambda_1-\lambda_2)-(\lambda_3-\lambda_2) = \lambda_1-\lambda_3$, we have the following algebraic relation,
$$
 \frac{{\bf N}_j}{K_{j}M_j}(r_j+ M_{j}q_j)-\frac{{\bf N}_k}{K_{k}M_k}(r_k+ M_{k}q_k)=\frac{{\bf N}_{\ell}}{K_{{\ell}}M_{\ell}}(r_\ell+ M_{{\ell}}q_\ell).
$$
It follows that
    $$
{\bf N}_jK_{k}K_{{\ell}}M_kM_{\ell}(r_j+ M_jq_j)-  {\bf N}_{k}K_{j}K_{{\ell}}M_{j}M_{\ell}(r_k+ M_kq_k)= {\bf N}_{\ell}K_{j}K_{{k}}M_{j}M_k(r_\ell+ M_{\ell}q_\ell).
    $$
    Hence,
   \begin{equation}\label{eq1.2}
\begin{aligned}
&  {\bf N}_jK_{k}K_{{\ell}}M_kM_{\ell}r_j-{\bf N}_{k}K_{j}K_{{\ell}}M_{j}M_{\ell}r_k-{\bf N}_{\ell}K_{j}K_{{k}}M_{j}M_kr_\ell  \\
=&M_jM_kM_{\ell}\cdot ( {\bf N}_{\ell}K_{j}K_{{k}}q_\ell-{\bf N}_{k}K_{j}K_{{\ell}}q_k-{\bf N}_jK_{k}K_{{\ell}}q_j).
\end{aligned}
\end{equation}

\medskip

In the first case, if all $M_j=2$, then all $0<r_i<2$ which means all $r_i=1$.  (\ref{eq1.2}) is reduced to
$$
 {\bf N}_jK_{k}K_{{\ell}}-{\bf N}_{k}K_{j}K_{{\ell}}-{\bf N}_{\ell}K_{j}K_{k} =  2\cdot ( {\bf N}_{\ell}K_{j}K_{{k}}q_\ell-{\bf N}_{k}K_{j}K_{{\ell}}q_k-{\bf N}_jK_{k}K_{{\ell}}q_j).
$$
The right hand side  is an even number. However, as all $K_{j}$ are odd numbers, and all $N_j = 2K_j+1$ are odd,, each term in the left hand side of (\ref{eq1.2}) must be odd and hence the left hand side is an odd number overall. This is a contradiction. Hence, we cannot have a mutually orthogonal set of cardinality large than $2$.

\medskip

\end{pf}

\medskip

We are now ready to construct the example in Theorem \ref{th0.4}.

\medskip

\noindent{\it Proof of Theorem \ref{th0.4}.} Let $p$ be a prime number of the form $4k+3$. It is well known that there are infinitely many primes of such form by the Dirichlet theorem. Writing
$$
p^j = 2K_j+1,
$$j
we claim that $K_j$ is an odd number whenever $j$ is odd. Expanding in binomial theorem, we have for some integer $L_j$,
$$
K_j = \frac{p^j-1}{2} = \frac{(4k+3)^j-1}{2} =\frac{4L_j+3^j-1}{2}.
$$
It suffices to show that $(3^j-1)/2$ is an odd number if $j$ is odd. But from Binomial expansion, $
\frac{3^j-1}{2} = \frac{(2+1)^j-1 }{2}= 2^{j-1}+j2^{j-2}+...\binom{j}{2}2+j.$ This shows that $\frac{3^j-1}{2}$ is an odd number.

\medskip

Letting $N_j = p^{2j-1} = 2K_{2j-1}+1$ where $p$ is a prime number of the form $4k+3$ and $B_{j} = \{0,K_{2j-1}\}$. From Example \ref{Example2.3}, we have an almost-Parseval-frame tower. By Theorem \ref{th0.2}(b),  the associated measure $\mu$ is frame spectral. On the other hand, it is non-spectral by Proposition \ref{prop1.2}. \qquad$\Box$

\medskip

We end this section with a remark.

\begin{Rem}
{\rm In \cite{[HLL]}, it was proved that if $\nu$ is a spectral measure on $[0,1]$ with spectrum inside ${\Bbb Z}$, then any $\mu = \nu\ast \delta_{A}$, with $A\subset{\mathbb Z}$, is a frame spectral measure and some of them are not spectral. In view of Theorem \ref{th0.1}, the measure $\mu$ we constructed cannot be of the form $\nu\ast\delta_{A}$, where  $\nu$ is spectral and $\delta_A$ is a discrete measure supported on some set $A$. If it was the case, then}
$$
\widehat{\mu}(\xi) = \widehat{\nu}(\xi)\widehat{\delta_A}(\xi)
$$
{\rm and this would have implied that any mutually orthogonal set of $\nu$ must be mutually orthogonal set of $\mu$ and hence cardinality of such sets for $\mu$ for be infinite, which is a contradiction by Proposition \ref{prop1.2} we just proved.}
\end{Rem}

\medskip

\section{Appendix: Hausdorff dimension}

In this section, we study the Hausdorff dimension, denoted by $\mbox{dim}_H$, of the support of $\mu$, which is an important question in fractal geometry. We refer the reader to \cite{[Fal]} the definition of Hausdorff dimension. Given a sequence of positive integers $M_j$ and a sequence of numbers  $r_j$. Suppose that $0<r_j<1$, $M_j\ge 2$ and $r_jM_j\le1$ for all $j$. For $k\ge 1$, we let $D_0 = \emptyset$, $D_k = \{(i_1,...,i_k): i_j\in \{0,1,...,M_j-1\}\}$. For ${\bf i}\in D_k$ and ${\bf j}\in D_{\ell}$, ${\bf i}{\bf j}\in D_{k+\ell}$ is the standard concatenation of two words. For each $\sigma\in \bigcup_{k=1}^{\infty}D_k$, we define an interval  $J_{\sigma}$. We say that
$$
E= \bigcap_{k=1}^{\infty}\bigcup_{\sigma\in D_k} J_{\sigma}
$$
is a {\it homogeneous Moran set} if the following conditions are satisfied.

\begin{enumerate}
\item $J_{\emptyset} = [0,1]$. For any $\sigma\in D_k$, $J_{\sigma0}$,...,$J_{\sigma(M_{k+1}-1)}$ are subinterval of $J_{\sigma}$ enumerated from left to right and $J_{\sigma i}\cap J_{\sigma j}$ has  intersects at most one point.
\item For any $k\ge 1$, for any $\sigma\in D_{k-1}$ and $j\in \{0,1,...,M_k-1\}$,
$$
r_k = \frac{|J_{\sigma j}|}{|J_{\sigma}|}.
$$
$|\cdot|$ denotes the Lebesgue measure of the interval.
\item For any $\sigma\in D_k$, the gaps between $J_{\sigma i}$ and $J_{\sigma (i+1)}$ are equal in length, the left endpoint of $J_{\sigma 0}$ is equal to the left endpoint of $J_{\sigma}$ and the right endpoint of $J_{\sigma (M_{k+1}-1)}$ is equal to the right endpoint of $J_{\sigma}$.
\end{enumerate}
It was shown that \cite{[FWW]} (see also \cite[Proposition 3.1]{[FLW]}) that the Hausdorff dimension of $E$ is equal to
$$
\mbox{dim}_H(E) = \liminf_{j\rightarrow\infty}\frac{\log (M_1...M_j)}{-\log(r_1...r_j)}.
$$
Turning to our case where $N_j = M_jK_j+\alpha_j$ and $B_j = \{0, K_j,...,K_j(M_j-1)\}$, the support of the measure $\mu$ is $$
K_{\mu} = \left\{\sum_{j=1}^{\infty} \frac{i_j K_j}{N_1...N_j}: i_jK_j\in B_j\right\} = \bigcap_{k=1}^{\infty}\bigcup_{\sigma\in D_k} J_{\sigma}
$$
where $J_{\sigma} = [ \sum_{j=1}^{k} i_jK_j(N_1...N_j)^{-1}, \sum_{j=1}^{k} i_jK_j(N_1...N_j)^{-1}+(N_1...N_k)^{-1}]$ and $\sigma = (i_1,...,i_k)$. Note that The support $K_{\mu}$ is contained in the interval $[0, \rho]$, where $\rho =\sum_{j=1}^{\infty} (K_j(M_j-1))(N_1...N_j)^{-1} . $ By a simple rescaling, $C = \rho^{-1} K_{\mu}$. It is easy to see that $C$ is actually a homogeneous Moran set with $r_1 = 1/N_1$ and $r_k = (N_1...N_{k})^{-1}/(N_1...N_{k-1})^{-1} = 1/N_k$. Hence, we have thus proved

\begin{Prop}
$$
\mbox{dim}_H(K_\mu) = \liminf_{j\rightarrow\infty}\frac{\log (M_1...M_j)}{\log(N_1...N_j)}.
$$
\end{Prop}

\medskip

We now compute the Hausdorff dimension of  some frame spectral measures $\mu$.

\medskip

\begin{Example}
{\rm In Example \ref{Example2.3}, $N_j = p^j = 2K_j+1$ and $M_j = 2$ for all $j$. In this case, the Hausdorff dimension}
$$
\mbox{dim}_H(K_\mu) = \liminf_{j\rightarrow\infty}\frac{\log 2^j}{\log(p^{1+2+...+j})} = \liminf_{j\rightarrow\infty}\frac{2j\log 2}{ j(j+1)\log p} =0 .
$$
{\rm The non-spectral measure with Fourier frame given in Theorem \ref{th0.4} is a special case of this type, and thus the support has Hausdorff dimension 0.}
  \end{Example}

This example shows that frame spectral measure can be very `` thin", similar situation happens for spectral measures \cite{[DaS]}. The following example shows that our construction does give frame spectral measures with positive Hausdorff dimension.

\begin{Example}
{\rm In Example \ref{Example2.4}, $M_j =K_j^{\beta}$ and $\alpha_j =1$ (i.e. $\gamma=0$) for all $j$. Then $N_j = K_j^{1+\beta}+1$. In this case, the Hausdorff dimension}
$$
\mbox{{\rm dim}}_H(K_\mu) = \liminf_{j\rightarrow\infty}\frac{\log ((K_1...K_j)^{\beta})}{\log((K_1^{1+\beta}+1)...(K_j^{1+\beta}+1))} .
$$
{\rm As $\lim_{x\rightarrow\infty}\log(1+x)/\log x =1$, for $x$ large, $C^{-1}\log x\leq \log(1+x)\leq C\log x$ for some  constant $C>0$. Hence, if $K_j$ is large enough,}
$$
C^{-1}\log((K_1...K_j)^{1+\beta})\leq\log((K_1^{1+\beta}+1)...(K_j^{1+\beta}+1))\leq C\log((K_1...K_j)^{1+\beta})
$$
{\rm This implies that }
$$
\begin{aligned}
\frac{\beta}{C(1+\beta)}=\liminf_{j\rightarrow\infty}\frac{\beta\log ((K_1...K_j))}{C(1+\beta)\log((K_1...K_j))} \leq&\mbox{{\rm dim}}_H(K_\mu) \\ \leq& \liminf_{j\rightarrow\infty}\frac{\beta\log ((K_1...K_j))}{C^{-1}(1+\beta)\log((K_1...K_j))}
\\=& \frac{\beta}{(1+\beta)C^{-1}}.\\
\end{aligned}
$$
{\rm Hence, the support of the frame spectral measure has Hausdorff dimension at least $\frac{\beta}{C(1+\beta)}>0$.}
  \end{Example}

\end{document}